\documentclass[10pt, oneside, a4paper]{article}

\usepackage{amsmath}
\usepackage{amssymb}
\input xy
\xyoption{all}

\newtheorem{theorem}{Theorem}
\newtheorem{lemma}[theorem]{Lemma}
\newtheorem{cor}[theorem]{Corollary}
\newtheorem{prop}[theorem]{Proposition}

\newtheorem{defn}[theorem]{Definition}

\newtheorem{rem}[theorem]{Remark}
\newtheorem{ex}[theorem]{Example}

\newcommand{\ac}{\`a }

\def\W{\Omega^{(n)}}

\newcommand{\la}{\lambda}
\newcommand{\Om}{\Omega}
\newcommand{\Omn}{\Omega^{(n)}}
\newcommand{\Omm}{\Omega^{(m)}}
\newcommand{\OmN}{[\Omega]^n}
\newcommand{\cvd}{\hspace*{\fill}
    {\rm \hbox{\vrule height 0.2 cm width 0.2cm}}}
\newcommand{\N}{{\mathbb N}}

\newcommand{\K}{{\mathbb K}}
\DeclareMathOperator{\Aut}{Aut}
\DeclareMathOperator{\Sym}{Sym}
\DeclareMathOperator{\Out}{Out}

\DeclareMathOperator{\Ker}{ker} \DeclareMathOperator{\Id}{id}

\title{Almost-free finite covers}

\author{Elisabetta Pastori \\
Dipartimento di Matematica, Universit\ac di Firenze, \\
Viale Morgagni 67A, 50134, Firenze, Italy. \\
pastori@math.unifi.it \\}

\begin{document}
\date{}

\maketitle

\begin{abstract}
Let $W$ be a first-order structure and $\rho$ be an $\Aut(W)$-congruence on $W$.
In this paper we define the {\emph{almost-free}} finite covers of $W$ with respect to $\rho$, and we show how to construct them. These are a generalization of free finite covers.\\
A consequence of a result of \cite{EH} is that any finite cover of $W$ with binding groups all equal to a simple non-abelian permutation group is  almost-free with respect to some $\rho$ on $W$. Our main result  gives a description (up to isomorphism) in terms of the $\Aut(W)$-congruences on $W$ of the kernels of principal finite covers of $W$ with bindings groups  equal at any point to a simple non-abelian regular permutation group $G$.  Then we analyze almost-free finite covers of  $\Omn$, the set of ordered $n$-tuples of distinct elements from a countable set $\Omega$, regarded as a structure with $\Aut (\Omn)=\Sym(\Om)$ and we show a result of biinterpretability.\\
The material here presented addresses a problem which arises in the context of  classification of totally categorical structures.
\end{abstract}

\section{Introduction}

Given an countable set $W$, consider the natural action of the
symmetric group $\Sym (W)$ on $W$. This action yields a topology on
$\Sym(W)$ in which pointwise stabilizers of finite sets give a base
of open neighborhoods of the identity. 
Let $\Upsilon$ be a closed
subgroup of $\Sym(W)$ that acts transitively on $W$ and $G$  a
finite group acting on a finite set $\Delta$. Consider the
projection $\pi:\Delta\times W\rightarrow W$ given by
$\pi(\delta,w)=w$. We denote by $G^W$ the set of all functions from
$W$ to $G$. Let $\mathcal{F}$ be the set of closed subgroups of
$\Sym(\Delta\times W)$ which preserve the partition of $\Delta\times
W$ given by the fibres of $\pi$. Every $F\in \mathcal{F}$ determines
naturally an induced map $\mu_F:F\rightarrow \Sym(W)$. Additionally
we require that, for all $F\in \mathcal{F}$,
$\mu_F(F)=\Upsilon$ and the permutation groups induced repectively
by $F$ and $\Ker\mu_F$ on $\pi^{-1}(w)$, for all $w\in W$, are both equal  to $G$.  Let $\mathcal{K}=\{\Ker\mu_F$,
$F\in \mathcal{F}\}$.  In this paper we will deal with the following\\

{\bf{Problem:}} Given $G$ and $\Upsilon$, find a description of 
the elements belonging to $\mathcal{K}$.\\

This problem, which is here formulated in terms of infinite
permutation groups, is motivated by questions arising in model
theory concerning finite covers (see \cite{EMI}).

\begin{defn}\label{finite cover}
Let $C$ and $W$ be two first-order structures. A finite to-one
surjection $\pi:C\rightarrow W$ is a {\emph{finite cover}} if its
fibres form an Aut$(C)$-invariant partition of $C$, and the induced
map $\mu:\textrm{Aut}(C)\rightarrow \Sym(W)$, defined by
$\mu(g)(w)=\pi(g\pi^{-1}(w))$, for all $g\in \textrm{Aut}(W)$ and
for all $w\in W$, has image Aut$(W)$.\\
We shall refer to the kernel of $\mu$ as the {\emph{kernel}} of the
finite cover $\pi$. If $\pi:C\rightarrow W$ is a finite cover, the
{\emph{fibre group}} $F(w)$ at $w\in W$ is the permutation group
induced by Aut$(C)$ on $\pi^{-1}(w)$ . The {\emph{binding group}}
$B(w)$ at $w\in W$ is the permutation group induced by the kernel on
$\pi^{-1}(w)$.
\end{defn}
Using the terminology of finite covers, the problem above can be
stated in the following equivalent version: given a finite group $G$
and a first-order structure $W$ with automorphism group $\Upsilon$,
describe the kernels of the finite covers of $W$ with $F(w)=B(w)=G$ at any point, which have $\Delta\times W$ as domain of the covering structures and are the projection on the second coordinate.\\
A more detailed commentary on finite covers and this problem is
given in the last section. However, we avoid the model-theoretic
methods using rather infinite permutation
groups techniques.\\

 In \cite{AZ2}
Ahlbrandt and Ziegler described the subgroups $K\in \mathcal{K}$, when $G$ is an abelian permutation group. In this case $G^W$, the group of the function from $W$ to $G$, is an $\Upsilon$-module with  $f^\upsilon(w)=f(\upsilon^{-1}w)$, where $\upsilon\in \Upsilon$ and $f\in G^W$ and the kernels in $\mathcal{K}$ are profinite
$\Upsilon$-modules. They proved that  $\mathcal{K}$ is exactly  the set of closed $\Upsilon$-submodules of $G^W$ .\\
In this paper, we deal with the case when $G$ is a simple {\emph{non
abelian}} regular  permutation group. Under this hypothesis our main result, which is
stated and proved in Section 3), gives a description of the elements
of $\mathcal{K}$ in terms of the $\Upsilon$-congruences on $W$. A
key ingredient in the proof is a result of Evans and
Hrushovski (\cite{EH}, Lemma 5.7).\\
Previous results are the following. In \cite{Zi}, Ziegler described the groups $K\in \mathcal{K}$ in the
case when $W$ is a countable set $\Omega$ and
$\Upsilon=\Sym(\Omega)$ (the disintegrated case), for any group $G$.
Increasing the complexity of the set $W$, it seems not possible to
give a general description of the groups $K\in \mathcal{K}$ not
depending on the group $G$. For example, if $W$ is the set on
$n$-subsets from a countable set $\Omega$, $\Upsilon=\Sym(\Omega)$
and $G$ is a cyclic group of order a prime $p$, then the groups
$K\in \mathcal{K}$ are an intersection of kernels of certain
$\Upsilon$-homomorphism, as it is described in \cite{Gr}. While if
$G$ is a simple non abelian group, then $\mathcal{K}=\{G,G^W\}$ (see
corollary \ref{primitivita}). In Section 4 we analyze the special
case in which given a countable set $\Omega$, $W$ is defined as the
subset of the $n$-fold cartesian product $\Omn$ whose elements are
$n$-tuples with pairwise distinct entries. Defining $\Upsilon$ as
$\Sym(\Omega)$, in Proposition 16 and 17  we give an explicit
description of the equivalence classes of the
$\Sym(\Omega)$-congruences on $\Omn$. In these Propositions  we see  that the blocks for $\Sym(\Om)$ in $\Omn$ can be either of finite or of infinite cardinalities. Proposition 23 shows that if $\pi:C\rightarrow \Omn$ is a cover of $\Omn$ with $\Aut(C)$ in $\mathcal{F}$ and $G$ equal to a simple-non abelian finite group such that the kernel of $\pi$ determines a $\Sym(\Om)$-congruence on $\Omn$ (in the sense of Lemma \ref{pregeometry}) with classes of finite cardinality, then, for every $m\in \N$ greater then $n$, there exists a finite cover $\pi':C'\rightarrow \Omm$ bi-interpretable with $\pi$ with binding groups and fibre groups both equal to $G$ at any point and kernel that determines a $\Sym(\Om)$-congruence on $\Omm$ with classes of infinite cardinality.\\
In section 5.3 we define the almost-free finite covers. A posteriori we see that the results of sections 3 and 4  concern examples of almost-free finite covers with binding groups equal to the fibre groups at any point.  Let $W$ be a transitive structure, $\rho$ be an $\Aut(W)$-congruence on $W$ and $[w_0]$ be a congruence class.\\ 
An almost-free finite cover $\pi$ of $W$ w.r.t $\rho$ is a finite cover whose permutation group induced by its kernel on the union of the fibres of $\pi$  over  $[w_0]$ is isomorphic to the binding group at $w_0$, while the permutation group induced on the fibres over two elements not in the same congruence class is the direct product of the two respective binding groups.  This definition generalizes the definition of free finite cover. More in detail a free finite cover of $W$ is an almost- free finite cover of $W$ with respect to the equality. In Proposition 26 we show how to construct an almost-free finite cover. The proof uses Lemma 2.1.2 of \cite{EMI}.

\section{General results}
\begin{defn}
A {\emph{pregeometry}} on a set $X$ is a relation between elements
$x\in X$ and finite subsets $X_0\subset X$, called
{\emph{dependence}}, which satisfies:
\begin{itemize}
  \item Reflexivity : $x$ is dependent on $\{x\}$;
  \item Extension: $x$ depends on $X_0$ and $X_0\subseteq X_1$
  implies $x$ depends on $X_1$;
  \item Transitivity: $x$ is dependent on $X_0$ and every $y\in X_0$
  is dependent on $X_1$ implies $x$ is dependent on $X_1$;
  \item Symmetry: $x$ is dependent on $X_0\cup \{y\}$ but not on
  $X_0$, implies $y$ is dependent on $X_0\cup \{x\}$.
\end{itemize}
\end{defn}
\begin{rem}
A classical example of a pregeometry is a vector space with linear
dependency.
\end{rem}
If $\Om$ is any set then there is a natural topology on
$\Sym(\Om)$ which makes it into a topological group. The open sets
are unions of cosets of pointwise stabilizers of finite subsets of
$\Om$. We then make any permutation group $P$ on $\Om$ into a
topological group by giving it the relative topology. If $\Om$ is
countable the topology is metrisable.\\ 
From now on $W$ stands for a countable set, $\Upsilon$ for a closed
subgroup of $\Sym(W)$ that acts transitively on $W$ and $G$ for a
finite group acting on a finite set $\Delta$. Consider the
projection $\pi:\Delta\times W\rightarrow W$ given by
$\pi(\delta,w)=w$. We denote by $G^W$ the set of all functions from
$W$ to $G$. Let $\mathcal{F}$ be the set of closed subgroups of
$\Sym(\Delta\times W)$ which preserve the partition of $\Delta\times
W$ given by the fibres of $\pi$. Every $F\in \mathcal{F}$ determines
naturally an induced map $\mu_F:F\rightarrow \Sym(W)$. Additionally
we require that, for all $F\in \mathcal{F}$,
$\mu_F(F)=\Upsilon$ and the permutation groups induced repectively
by $F$ and $\Ker\mu_F$ on $\pi^{-1}(w)$, for all $w\in W$, are both equal  to $G$. We notice that the wreath product $G Wr_W \Upsilon$ in its imprimitive action on $\Delta\times W$ belongs to $\mathcal{F}$.\\
 It is easy to see that, with the above topology, $G^W$ is a compact subgroup of $\Sym(\Delta\times W)$ and $\Ker\mu_F$ are
closed subgroups of $G^W$ and that $\mu_F$ are
continuous and open maps (Lemma 1.4.2, \cite{EMI}). We introduce now a notion of isomorphism
among the elements of $\mathcal{F}$. We say that $F_1$ and $F_2$ are
$\emph{isomorphic}$ if there exists a bijection $\phi:\Delta \times
W\rightarrow \Delta \times W$ which sends
$\phi(\pi^{-1}(w))=\pi^{-1}(w)$, for all $w\in W$ and such that the
induced map $f_\phi:\Sym(\Delta\times W)\rightarrow \Sym(\Delta
\times W)$ sends $F_1$ to $F_2$. Let $\mathcal{K}=\{\Ker\mu_F$,
$F\in \mathcal{F}\}$.  We now introduce the following equivalence
relation $R$ on $\mathcal{K}$: $\Ker\mu_{F_1}R\Ker\mu_{F_2}$ if and
only if $F_1$ is isomorphic to $F_2$ and we denote the
$R$-equivalence class of an arbitrary $K\in \mathcal{K}$ by $[K]$.
(We shall say that $\Ker\mu_{F_1}$ is $\emph{isomorphic}$ to
$\Ker\mu_{F_2}$ if $F_1$ is isomorphic to
$F_2$.)\\

Take $K\in \mathcal{K}$ and $w_1,\dots, w_k\in W$. We define
$$K(w_1,\dots,w_k)=\{f|_{\{w_1,\dots,w_k\}}\, | \, f\in K
\}$$ and, for simplicity, we shall refer to $K(w_1,\dots,w_k)$ as
$K$ restricted to $w_1,\dots,w_k$.
\begin{defn}
Suppose $w_1,\dots,w_k,w$ belong to $W$.We say that $w$ depends on
$w_1\dots,w_k$ and write $w\in \emph{cl}(w_1,\dots,w_k)$, if
$$
K(w,w_1,\dots,w_k)\cong G
$$
\end{defn}
\begin{lemma}[\cite{EH}, Lemma 5.7]\label{pregeometry}
Let $K\in \mathcal{K}$ and $w_1,\dots,w_k,w\in W$. Then
$(W,\emph{cl})$ is a $\Upsilon$-invariant pregeometry. If $G$ is a
simple non-abelian finite group, then $(W,\emph{cl})$ reduces to
an equivalence relation.
\end{lemma}
The lemma states that, if $G$ is non-abelian and $w$ depends on
$w_1,\dots,w_k$, then there is an $i\in \{1,\dots,k\}$ such that
$w$ depends on $w_i$ and $(W,\textrm{cl})$ is a
$\Upsilon$-congruence.
\begin{cor}\label{primitivita}
If $\Upsilon$ acts primitively on $W$ and $G$ is a simple
non-abelian finite group, then $\mathcal{K}=\{G,G^W\}$.
\end{cor}
Here there are some results on topological groups that will be
useful in the next section.

\begin{lemma}\label{trivialita}
Let $G$ be a permutation group on an infinite set $\Omega$ with
the usual topology. A subgroup $H$ of $G$ is open in this topology
if and only if $H$ contains $G_{(\Gamma)}$ for some finite
$\Gamma$.
\end{lemma}

Take a typical basic open set of Aut$(C_1)=G_1$:
$$
(G_1)_{(F)}=\{g_1\in G_1:\,g_1(f)=f,\,\forall f\in F\}
$$
for some finite $F\subset C_1$. Let $a_i\in \Om$, for $i\in
\{1,\dots,n\}$, and $\la\in \Lambda$. The preimage of
$(G_1)_{(F)}$ under $\Phi$ is $ \Phi^{-1}((G_1)_{(F)})=\{g\in
\textrm{Aut}(C):\,g(f')=f',\, \forall f'\in F'\},$ where
$F'=\{(\la,a_1,\dots,a_m):(\la,a_1,\dots,a_n)\in F\}$ is a finite
set. Hence, $\Phi^{-1}((G_1)_{(F)})=\textrm{Aut}(C/F')$, which is
open and $\Phi$ is continuous.\\

\begin{prop}\label{compact}
Let $G$ be a topological group and let $H$ be a subgroup of $G$.
Then, if $G$ is compact and $H$ closed, $H$ is compact.
\end{prop}

For a proof of the previous proposition see for instance
\cite{Hig} Chapter 2, paragraph 8,10.
\begin{prop}\label{closedcompact}
Let $G$ be a topological group. Suppose $G$ is metrisable. Let $A$
be a compact subgroup of $G$ and $B$ a closed subgroup of $G$.
Then $AB$ and $BA$ are closed sets.
\end{prop}
\emph{Proof.} It is sufficient to show that $AB$ is closed. Let
$\{c_n\}_{n\in \mathbb{N}}$ be a sequence of elements of $AB$
which converges to $c$. We have $c_n=a_n b_n$, where $a_n\in A$
and $b_n\in B$. Since $A$ is compact, we can select from the
sequence $\{a_n\}_{n\in \mathbb{N}}$ a subsequence $\{a_{n_k}\}$
which converges to an element $a\in A$. We conclude from the
convergence of the sequences $\{c_{n_k}\}$ and $\{a_{n_k}\}$ that
the sequence $\{b_{n_k}\}$ converges to the element $a^{-1}c$,
which belongs to $B$, since $B$ is closed. Hence $c=a(a^{-1}c)\in
AB$ and the closure of the set $AB$ is established.\cvd

\section{Main Theorem}

 We will denote by $\mathcal{C}$  the set of all
$\Upsilon$-congruences on $W$. 
\begin{defn} Let $\rho\in
\mathcal{C}$. We define the subgroup of $G^W$
$$K_{\rho}=\{f:W \rightarrow G : \hbox{$f$ constant on $Y$, $\forall \, Y\in W/\rho$}\}.$$
\end{defn}

\begin{theorem} \label{congruenze e Kernel}
 Let $G$ be a simple non-abelian finite permutation group acting regularly on a finite set $\Lambda$. Then there exists a bijection $\Psi$ between
$\mathcal{C}$ and $\mathcal{K}/R$ given by $\Psi(\rho)=[K_{\rho}]$.
The inverse mapping $\Phi$ of $\Psi$ is given by $\Phi([K])=\rho_K$,
where $\rho_K$ is defined by:
$$
w_i\rho_K w_j\Leftrightarrow K(w_i,w_j)\cong G.
$$
\end{theorem}
{\emph{Proof. }}We first show that $\Psi$ maps $\mathcal{C}$
into $\mathcal{K}/R$.\\
Let $\rho\in \mathcal{C}$. Then $K_{\rho}$, is a subgroup of $G^W$.
First of all we embed $K_{\rho}$ into $G^W\rtimes \Upsilon$ in the
natural way:
\begin{displaymath}
\begin{array}{ccc}
K_{\rho} & \hookrightarrow & G^W \rtimes \Upsilon\\
f & \mapsto & (f, 1)
\end{array}
\end{displaymath}
and then we notice that $K_{\rho}$ is normalized by $\Upsilon$.
Indeed, given $\sigma\in \Upsilon$, we have that
$$
(\sigma(f),1)(\la,w):=(1,\sigma)\,(f,1)\,(1,\sigma^{-1})(\la,w)=(f(\sigma^{-1}
w)\la,w).
$$
Since $f\in K_{\rho}$, for every $w_i\in [w_j]_{\rho}$ in $W$ we
have $f(w_i)=f(w_j)$, but, since $\rho$ is a $\Upsilon$-congruence
on $W$, we have $f(\sigma^{-1} w_i)=f(\sigma^{-1} w_j)$, for
every $w_i\in [w_j]_{\rho}$ and so $(\sigma(f),1)\in K_\rho$.\\
Since $K_{\rho}$ is normalized by $\Upsilon$, we can consider the
group:
$$
H:=K_{\rho}\rtimes \Upsilon.
$$
This is a subgroup of $G^W\rtimes \Upsilon$ and if  $\mu: GWr_W
\Upsilon \mapsto \Upsilon $ is the  map defined by
$\mu(f,\gamma)=\gamma$,  we then have that $\mu(H)=\Upsilon$ and
$\Ker\mu=K_{\rho}$. In order to prove that $K_{\rho}$ is an element
of $\mathcal{K}$ it is sufficient to show that $H$ is a closed
subgroup of $G^W\rtimes \Upsilon$. Indeed, $G^W\rtimes \Upsilon$ is
closed in $\Sym(\Delta\times W)$.\\
The first step is to prove that $K_{\rho}$ is closed.\\
 The finite group $G$ has
 the discrete topology, while $G^W$ has the product topology.  An element  $f\in G^W$ is a function from
$W$ to $G$.
 The $w$-projection map is the map
  $\pi_{w}:G^W\rightarrow G$ such
 that $\pi_{w}(f)=f(w)$. A basis for the product
 topology on $G^W$ is the family of all finite intersections of
 $\pi_{w}^{-1}(U)$, where $U$ is an open subset of $G$. In this topology the maps $\pi_{w}$ are continuous. Hence,
 a
 member of this basis is of the form
 $$
 \bigcap\{\pi_{w}^{-1}(U_{w})\,:\,w\in F\}
 $$
 where $F$ is a finite subset of $W$.\\
Let $[w]_{\rho}$ be a $\rho$-class and $g$ an element of the simple
finite group $G$. By the  continuity of $\pi_w$, $\pi^{-1}_w(g)$ is
a closed subset of $G^W$ . Let
$$M_{[w]_{\rho}}(g):=\bigcap_{v\in
[w]_{\rho}}\pi^{-1}_{v}(g).$$ Then $M_{[w]_{\rho}}(g)$ is a closed
set in $G^W$. We consider next

$$
\bigcup_{g\in G}M_{[w]_{\rho}}(g)
$$
and this is still a closed subset of $K_0$. Then, if $\Sigma$ is the
set of all the equivalence classes of $\rho$,
$$
K_{\rho}=\bigcap_{[w]_{\rho}\in \Sigma}\bigcup_{g\in
G}M_{[w]_{\rho}}(g)
$$
and so $K_{\rho}$ is closed in $K_0$.\\
Since $K_\rho$ is a closed subgroup of the compact group $G^W$,
$K_\rho$ is compact by Proposition \ref{compact}. By Proposition
\ref{closedcompact}, $H=K_\rho \rtimes \Upsilon$ is closed. Thus, we
have shown that $\Psi$ maps $\mathcal{C}$ to $\mathcal{K}/R$.\\

It's easy to see that the map $\Phi$ is well defined. Finally, Lemma
\ref{pregeometry} shows that $\Phi([K])\in \mathcal{C}$.\\
In order to prove that $\Psi$ is a bijection, we  show  that $\Phi
\circ \Psi=\Id$
on $\mathcal{C}$.\\
Let $\rho$ be a $\Upsilon$-congruence on $W$ and let
$\Phi([K_{\rho}])=\bar{\rho}$. We want
to prove that $\rho=\bar{\rho}$.\\
Let $w_i,w_j\in W$ such that $w_i\, \rho \, w_j$, then for every
$f\in K_\rho$, $f$ is constant on the equivalence class
$[w_i]_\rho$, i.e. $f(w_i)=f(w_j)$. Hence, $K_\rho(w_i,w_j)\cong G$
and $[w_i]_\rho \subseteq [w_i]_{\bar{\rho}}$. Vice versa, let
$w_i\in W$ and suppose there exists $w_j\in W$ such that $w_j\notin
[w_i]_\rho$, but $w_j\in [w_i]_{\bar{\rho}}$. Since $w_j\notin
[w_i]_\rho$,   there exists an $f\in K_\rho$ such that $f(w_i)=g$
and $f(w_j)=1$, where $g\in G$ and $g\neq 1$. Then
$K_\rho(w_i,w_j)= G\times G$ and this yields a contradiction.\\
We shall finally prove that  $\Psi \circ \Phi=\Id$.\\
Let $K\in \mathcal{K}$ , $\Phi([K])=\rho_K$ and
$$\Psi(\Phi([K]))=[K_{\rho_K}].$$

 Let  $w_j\in [w_i]_{\rho_K}$. Since
$K(w_i,w_j)\cong G$, it means that there exist automorphisms
$\alpha_{w_i}, \alpha_{w_j}\in \Aut(G)$ such that, for every $f\in
K$, there exists $g\in G$ such that $f(w_i)=\alpha_{w_i}(g)$ and $f(w_j)=\alpha_{w_j}(g)$. We denote by $N_{\Sym(\Delta)}(G)$ the normalizer of $G$ in
$\Sym(\Delta)$.   Since $G$ acts regularly on $\Delta$,  for every $w\in W$
there exists $n_w$ belonging to  $N_{\Sym(\Delta)}(G)$ such that
$\alpha_{w}(g)=n_w^{-1}gn_w$, for $g\in G$. Consider the function
$n:W\rightarrow N_{\Sym(\Delta)}(G)$ given by $n(w)=n_w$. Let $F_{\rho_K}\in \mathcal{F}$ be a closed subgroup of $\Sym(\Delta\times W)$ such that $K_{\rho_K}=F_{\rho_K}\cap G^W$.  Since $F_{\rho_K}$ is closed,  $n^{-1}F_{\rho_K}n$ is closed. In fact,  $n^{-1}F_{\rho_K}n\in \mathcal{F}$ and
$$K=n^{-1}K_{\rho_K}n=n^{-1}F_{\rho_K}n\cap G^W.$$
Since $n$ is a bijection of $\Delta \times W$ which preserves the fibres of $\pi$, we have that  $n^{-1}F_{\rho_K}n$ is isomorphic to  $F_{\rho_K}$ and then $[K]=[K_{\rho_K}].$
\cvd
\begin{rem}\label{rem1}
It is clear by the previous proof that in every class $[K]\in \mathcal{K}/R$ there exists $\bar{K}\in [K]$ which is constant on the equivalence classes of $\Phi([K])$.
\end{rem}

\section{Special case}
Let $H$ be a group acting on a set $X$, $a\in X$ and
$\Delta\subseteq X$. We denote by $ a^H=\{ha:\, h\in H\}$,
 by $ H_{(\Delta)}$ the pointwise stabilizer of $\Delta$ in $H$ and by $H_{\{\Delta\}}$
 the setwise stabilizer of $\Delta$ in $H$.  We recall the
following theorem, whose proof can be found in \cite{DiMo}.
\begin{theorem}[\cite{DiMo},
Theorem 1.5A] \label{theorem:box} Let $G$ be a group which acts
transitively on a set $\Omega$, and let $\alpha\in \Omega$. Let
$\mathcal{D}$ be the set of blocks $\Delta$ for $G$ containing
$\alpha$, let $\mathcal{H}$ denote the set of all subgroups $H$ of
$G$ with $\textrm{G}_{\alpha}\leq H$. There is a bijection $\Psi$
from $\mathcal{D}$ onto $\mathcal{H}$ given by
$\Psi(\Delta):=\textrm{G}_{\{\Delta\}}$ whose inverse mapping
$\Phi$ is given by $\Phi(H):=\alpha^H$. The mapping $\Psi$ is
order preserving in the sense that if $\Delta,\Theta\in
\mathcal{D}$ then $\Delta\subseteq \Theta\Longleftrightarrow
\Psi(\Delta)\leq \Psi(\Theta)$.
\end{theorem}
From now on let $W$ be $\Omega^{(n)}$, the set of ordered $n$-tuples
of distinct elements of the countable set $\Omega$. Let
$\Upsilon=\textrm{Sym}(\Omega)$ act on $\Omn$ in the natural way:
let $\sigma\in \Sym(\Om)$, then
$\sigma(a_1,\dots,a_n)=(\sigma(a_1),\dots,\sigma(a_n))$. In the
sequel  we denote $\textrm{Sym}(\Omega)$ by $S$ when
$\textrm{Sym}(\Omega)$ acts on $\Omega$.  Let $\rho$ be a
$\Upsilon$-congruence, and $\Delta\subseteq \Omega^{(n)}$ be the
equivalence class of $\rho$ containing the element
$\alpha=(a_1,\dots,a_n)$. We will refer to $\Delta$ as a block of
imprimitivity containing $\alpha$.
\begin{defn}
Let $\alpha=(a_1,\dots,a_n)\in \Omn$. We define
$$
\emph{supp}(\alpha):=\{a_1,\dots,a_n\}.
$$
\end{defn}
By Theorem $\ref{theorem:box}$, the subgroup
$\Upsilon_{\{\Delta\}}=\{x\in \Upsilon|\quad x\Delta=\Delta\}$
contains the stabilizer $\Upsilon_{\alpha}=S_{(a_1,\dots,a_n)}$. A
proof of the following lemma can be found in \cite{DiMo}.
\begin{lemma}[\cite{DiMo} Lemma 8.4B]\label{intersezione}
Let $\Sigma_1$ and $\Sigma_2$ be subsets of an arbitrary set
$\Omega$ such that $|\Sigma_1\cap \Sigma_2|=|\Sigma_1|\leq
|\Sigma_2|$. Then
$$
\langle
\textrm{Sym}(\Sigma_1),\textrm{Sym}(\Sigma_2)\rangle=\textrm{Sym}(\Sigma_1\cup
\Sigma_2),
$$
(we identify $\textrm{Sym}(\Sigma)$ with the pointwise stabilizer
of $\Omega\smallsetminus \Sigma$ ).
\end{lemma}
\begin{prop} \label{inclusioni}
Let $\alpha=(a_1,\dots,a_n)\in \Omn$. Let $\Delta \neq \Omn$ be a
block containing $\alpha$. Let $\{\Gamma_i\}_{i\in I}$ be the set
of finite subsets of $\Omega$
 such that
$$
\Upsilon_{\alpha}\leq S_{(\Gamma_i)}\leq \Upsilon_{\{\Delta\}}.
$$
Let $\Gamma=\bigcap_{i\in I}\Gamma_i$. Then
$$
\Upsilon_{\alpha}\leq S_{(\Gamma)}\leq \Upsilon_{\{\Delta\}}\leq
S_{\{\Gamma\}}.
$$
Moreover $\Gamma$ is finite and   $\Gamma \subseteq
\{a_1,\dots,a_n\}$.
\end{prop}
{\emph{Proof. }} We notice that the index set $I$ is non-empty:
for instance the set $\{a_1,\dots,a_n\}$ belongs to
$\{\Gamma_i\}_{i\in I}$. Moreover, it is finite since every $\Gamma_i\leq\{a_1,\dots,a_n\}$. In order to prove that
$\Upsilon_{\alpha}\leq S_{(\Gamma)}$ it is sufficient to notice
that for every $i\in I$, $\Gamma\subseteq \Gamma_i$. Then
$\Upsilon_{\alpha}\leq
S_{(\Gamma_i)}\leq S_{(\Gamma)}$, for every $i\in I$.\\
We use Lemma \ref{intersezione} to prove the inclusion
$S_{(\Gamma)}\leq \Upsilon_{\{\Delta\}}$. Let $\Sigma_i=\Omega
\smallsetminus \Gamma_i$, for $i\in I$. Then by Lemma
\ref{intersezione} we have $ \langle S_{(\Gamma_i)}, i\in I
\rangle=S_{(\bigcap_{i\in I}\Gamma_i)} $
and so $S_{(\Gamma)}\leq \Upsilon_{\{\Delta\}}$.\\
Notice that $\Gamma$ is the smallest subset of $\Omega$ such that
$\Upsilon_{\alpha}\leq S_{(\Gamma)}\leq \Upsilon_{\{\Delta\}}.$ We
want to prove the set $\Gamma$ has the smallest cardinality among
the finite sets $X$ of $\Omega$ such that $S_{(X)}\leq
\Upsilon_{\{\Delta\}}.$ Suppose not, then there exists a finite
subset of $\Omega$, say $\Sigma$, with $|\Sigma|\lneq |\Gamma|$
and $S_{(\Sigma)}\leq
\Upsilon_{\{\Delta\}}$.\\
If $\Gamma\cap \Sigma \neq \emptyset$, then by Lemma
\ref{intersezione}, we have
$$
\Upsilon_{\alpha}\leq S_{(\Gamma)} \leq S_{(\Gamma \cap \Sigma)}
\leq \Upsilon_{\{\Delta\}}
$$
and, since $\Gamma$ is the smallest subset of $\Omega$ such that
$\Upsilon_{\alpha}\leq S_{(\Gamma)}\leq \Upsilon_{\{\Delta\}},$ this yields a contradiction.\\
If $\Gamma \cap \Sigma= \emptyset$, then
$$\langle S_{(\Gamma)},S_{(\Sigma)}\rangle= S_{(\Gamma \cap
\Sigma)}=S \leq \Upsilon_{\{\Delta\}}
$$
but $\Upsilon_{\{\Delta\}} \neq S$, a contradiction. Thus, the set  $\Gamma$ has
the smallest cardinality  among the finite subsets $X$ of $\Omega$
such that $S_{(X)}\leq \Upsilon_{\{\Delta\}}.$\\
Let $x\in \Upsilon_{\{\Delta\}}$, then we have $
S_{(x\Gamma)}=x^{-1}S_{(\Gamma)}x \leq \Upsilon_{\{\Delta\}},$ and
so, applying again Lemma \ref{intersezione} we get  that
$\Upsilon_{\{\Delta\}}\geq \langle
S_{(\Gamma)},S_{(x\Gamma)}\rangle=S_{(\Gamma \cap x\Gamma)}$.
Thus, for all $x\in \Upsilon_{\{\Delta\}}$, $\Gamma=x\Gamma$ by
the minimality of $\Gamma$ and
$\Upsilon_{\{\Delta\}} \leq S_{\{\Gamma\}}$.\\
To prove that $\Gamma\subseteq \{a_1,\dots,a_n\}$ it is sufficient
to note that $S_{(\Gamma)}\geq \Upsilon_{\alpha}$, and the claim
follows.\cvd

As the following result shows, a $\rho$-class in $\W$ can be a
finite subset or an infinite subset of $\W$.
\begin{prop} \label{cardinality e stabilizzatori}
Let $\Delta\neq \W$ be the equivalence class of a
$\Upsilon$-congruence $\rho$
containing the element $(a_1,\dots,a_n)\in \Omn$. Then\\
{\emph{a) }}  $\Delta$ is finite if and only if
$S_{(a_1,\dots,a_n)}\leq \Upsilon_{\{\Delta\}} \leq
S_{\{a_1,\dots,a_n\}}$;\\
 \emph{b) }  $\Delta$ is a countably
 infinite set if and only if $S_{(\Gamma)}\leq \Upsilon_{\{\Delta\}} \leq S_{\{\Gamma\}}$, for some
finite set $\Gamma\varsubsetneqq  \{a_1,\dots,a_n\}$.

\end{prop}
{\emph{Proof. }}

  {\emph{a)}} Suppose $\Delta$ is a finite set in $\Omega^{(n)}$. If it doesn't exist any $\Gamma\varsubsetneqq \{a_1,\dots,a_n\}$ such that $S_{(a_1,\dots,a_n)}<S_{(\Gamma)}\leq \Upsilon_{\{\Delta\}}$, by Proposition \ref{inclusioni}, since $S_{(a_1,\dots,a_n)}\leq \Upsilon_{\{\Delta\}}$,  we have $S_{(a_1,\dots,a_n)}\leq \Upsilon_{\{\Delta\}} \leq
S_{\{a_1,\dots,a_n\}}$.\\
  Hence, suppose that
 there exists a finite set $\Gamma\varsubsetneqq \{a_1,\dots,a_n\}$ such that
$S_{(a_1,\dots,a_n)}\lneq S_{(\Gamma)}\leq \Upsilon_{\{\Delta\}}.$
Let $x\in S_{(\Gamma)}\leq \Upsilon_{\{\Delta\}}$, then $x
\Delta=\Delta$. Take $a_i\in \{a_1,\dots,a_n\}\setminus \Gamma$.
Then pick $a\in \Omega$ such that $a\notin \textrm{supp}(\delta)$,
for every $\delta\in \Delta$. By $k$-transitivity of $S$, for any
$k\in \mathbb{N}$, it is possible to choose an element $x$ in
$S_{(\Gamma)}$, such that $x(a_i)=a$. Then
$$x(a_1,\dots,a_i,\dots,a_n)=(x(a_1),\dots,a,\dots,x(a_n))\in
\Delta.$$ But this yields a contradiction, since $a\notin
\textrm{supp}(\delta)$, for every $\delta\in \Delta$.\\
In the other direction, if $S_{(a_1,\dots,a_n)}\leq \Upsilon_{\{\Delta\}} \leq
S_{\{a_1,\dots,a_n\}} $ then
$\Delta=(a_1,\dots,a_n)^{\Upsilon_{\{\Delta\}}}\subseteq
(a_1,\dots,a_n )^{S_{\{a_1,\dots,a_n\}}}$, and $| (a_1,\dots,a_n
)^{S_{\{a_1,\dots,a_n\}}}|$ is finite.\\
 {\emph{b) }} We now assume $\Delta$ is a countably infinite set. Suppose there does not exist
 any finite set $\Gamma\varsubsetneqq \{a_1,\dots,a_n\}$ such that
 $S_{(\Gamma)}\leq \Upsilon_{\{\Delta\}}$. By Theorem \ref{theorem:box} we have that $S_{(a_1,\dots,
 a_n)}\leq \Upsilon_{\{\Delta\}}$. Since for every  finite set $\Gamma\varsubsetneqq \{a_1,\dots,a_n\}$
 we have
 $S_{(\Gamma)}\nleq \Upsilon_{\{\Delta\}}$, then $\{a_1,\dots,a_n\}$ is
 the smallest subset of $\Omega$ such that $S_{(a_1,\dots,
 a_n)}\leq \Upsilon_{\{\Delta\}}$ and so, by Proposition \ref{inclusioni},
 $\Upsilon_{\{\Delta\}}\leq S_{\{a_1,\dots,a_n\}}$.
 Take an element $(b_1,\dots,b_n)$  of $\Delta$,
 such that $\{b_1,\dots,b_n\}\neq \{a_1,\dots,a_n\}$; as $\Delta$
 is infinite, this element there exists.
 By the $n$-transitivity of $S$, there exists an element $x\in S$ such that
 $x(a_1)=b_1,\dots,x(a_n)=b_n$.
 Then $x(a_1,\dots,a_n) \in \Delta$
 and so we have an element $x\in \Upsilon_{\{\Delta\}}$ but not in $S_{\{a_1,\dots,a_n\}}$.
 This yields a contradiction.\\
Conversely suppose    $\Gamma\varsubsetneqq \{a_1,\dots,a_n\}$,
and $S_{(\Gamma)}\leq \Upsilon_{\{\Delta\}} \leq S_{\{\Gamma\}}$.
Then $(a_1,\dots,a_n)^{S_{(\Gamma)}}\subseteq \Delta$, and since
$(a_1,\dots,a_n)^{S_{(\Gamma)}}$ is infinite, then $\Delta$ is
infinite.\cvd

\begin{rem}\label{groups}
If $|\Gamma|=n$, $n\geq 1$, then $S_{(\Gamma)}\unlhd
S_{\{\Gamma\}}$ and $S_{\{\Gamma\}}/S_{(\Gamma)}\cong Sym_n$ the
symmetric group on $n$ points. Given an element
$\alpha=(a_1,\dots,a_n)\in \Omn$ and a finite block $\Delta$
containing it, we have that $H=\Upsilon_{\{\Delta\}}$ satisfies
the following inclusions: $ S_{(\Gamma)}\leq H \leq
S_{\{\Gamma\}}\leq S$, where $\Gamma=\{a_1,\dots,a_n\}$. Then
$H/S_{(\Gamma)}$ is isomorphic to a subgroup of Sym$_n$. There
exists a bijection $\Theta$ between the subgroups of $\Sym_n$ and the
subgroups of $S_{\{\Gamma\}}$ which contain $S_{(\Gamma)}$.
\end{rem}

We shall denote by
$$\mathcal{K}_F=\{K\in \mathcal{K}|\, \hbox{$\rho_K$ has finite equivalence classes}\}.$$
\begin{prop}
Let $\mathcal{L}$ be the set of subgroups of Sym$_n$. Then there
exists a bijection
$$
\zeta:\mathcal{K}_F/R\rightarrow \mathcal{L}.
$$
\end{prop}
{\emph{Proof. }} By Theorem \ref{congruenze e Kernel}, it is
sufficient to find a bijection between the set of finite blocks
containing an element $\alpha=(a_1,\dots,a_n)$ and $\mathcal{L}$.
Let $\Delta$ be a finite block in $\Omn$ containing $\alpha$. We
have that
$$
\Upsilon_\alpha = S_{(\Gamma)}\leq \Upsilon_{\{\Delta\}} \leq
S_{\{\Gamma\}}\leq S
$$
where $\Gamma=$supp$(\alpha)$. Then by
Remark \ref{groups},  $\Upsilon_{\{\Delta\}}$ is the image by
$\Theta$ of a subgroup of Sym$_n$. If $\Delta_1\neq \Delta_2$ then
$\Upsilon_{\{\Delta_1\}}\neq \Upsilon_{\{\Delta_2\}}$. By 
Remark \ref{groups}, it follows that the map $\zeta$ is injective.
In the other direction, let $H\in \mathcal{L}$. By the remark \ref{groups},
$\Theta(H)$ is a subgroup $L$ of $S_{\{\Gamma\}}$ which contains
$\Upsilon_\alpha = S_{(\Gamma)}$. Then, by Theorem
\ref{theorem:box}, we have  a finite block $\alpha^L$ containing
$\alpha$.\cvd

\begin{prop}\label{blocco finito}
Let $\alpha=(a_1,\dots,a_n)\in \Omn$ and let
$\mathcal{D}_F^{\alpha}$ be the set of the finite blocks in $\Omn$
containing $\alpha$. Then the elements of $\mathcal{D}_F^{\alpha}$
are exactly the sets $\alpha^H$, where $H$ is a subgroup of
$\Sym\{a_1,\cdots,a_n\}$.
\end{prop}
{\emph{Proof.}}  Let $\Delta\in \mathcal{D}_F^{\alpha}$. Let $H'$
be the subgroup of Sym$(\Om)$ such that $\alpha^{H'}=\Delta$. Then
$$ \Upsilon_\alpha\leq H' \leq S_{\{\Gamma\}},
$$
where $\Gamma=\textrm{supp}(\alpha)$. Since
$S_{\{\Gamma\}}/S_{(\Gamma)}\cong Sym_n$ we have that $H'=H\times
\Sym(\Om\setminus \, \Gamma)$, where  $H$ is a  subgroup of
Sym$\{a_1,\cdots,a_n\}$. Then $\Delta=\alpha^H$.\\
Viceversa, taken a subgroup $H\leq \Sym\{a_1,\cdots,a_n\}$,
$\alpha^H=\alpha^{H\times \Sym(\Om \setminus \Gamma)}$. By Theorem
\ref{theorem:box} $\alpha^H$ is a block in $\Omn$ .\cvd

The same argument works for the following:
\begin{prop}
Let $\alpha=(a_1,\dots,a_n)\in \Omn$ and let
$\mathcal{D}_I^{\alpha}$ be the set of non-trivial infinite blocks
in $\Omn$ containing $\alpha$. Then the elements of
$\mathcal{D}_I^{\alpha}$ are exactly the sets $\alpha^{L\times
\Sym(\Om\setminus \Xi)}$, where $\Xi\subsetneq
\{\alpha_1,\cdots,\alpha_n\}$  and $L$ is a subgroup of
$\Sym(\Xi)$.
\end{prop}

Let us mention a little remark about Proposition \ref{blocco finito}.
Let $\alpha=(a_1,\dots,a_n)$. Denote $\Sym\{a_1,\dots,a_n\}$ by $\Sym_n$. Consider
the set
$$\alpha^{\Sym_n}=\{ \sigma(a_1,\dots,a_n), \sigma\in
\Sym_n\}.$$
Let $[K]\in \mathcal{K}_F/R$, and $\bar{K}\in [K]$ be the subgroup of $G^W$ such that is constant on the equivalence classes of $\Phi(K)$ (remind Remark \ref{rem1}). 
 By Proposition \ref{blocco finito} there
exists a subgroup $T$ of $\Sym_n$ such that  $K$ restricted to
$\Delta=\alpha^T$ is constant on it. The system of blocks
containing $\Delta$ is the set  $\{g\Delta, g\in \Sym(\Om)\}$. We look at the restriction
of $\bar{K}$ to the set $\alpha^{\Sym_n}$. This is the subgroup of
$G^{\alpha^{\Sym_n}}$ of the function from $\alpha^{\Sym_n}$ to
$G$ constant on the subsets $bT(\alpha)$, where $bT$ are the left
cosets of $T$ in $\Sym_n$.  We notice that the cardinalities of the finite
blocks in $\Omn$ are exactly the cardinalities of the subgroups of
Sym$_n$.

\section{Commentary}

\subsection{Finite Covers}

As is well known, a subgroup of $\Sym(W)$ is closed if and only if
it is the group of automorphisms of some first-order structure
with domain $W$ (see for instance Proposition (2.6) in \cite{Ca}). Thus we state the following definition.\\
A {\it{permutation structure}} is a pair $\langle W, G \rangle$,
where $W$ is a non-empty set (the {\it{domain}}), and $G$ is a
closed subgroup of Sym$(W)$. We refer to $G$ as the automorphism
group of $W$. If $A$ and $B$ are subsets of $W$ (or more generally
of some set on which Aut$(W)$ acts), we shall refer to Aut$(A/B)$
as the group of permutations of $A$ which extend to
elements of Aut$(W)$ fixing every element of $B$ and to $\Aut(A/\{B\})$ as the group of permutations of $A$ which extend to
elements of Aut$(W)$ stabilizing setwise the set $B$.\\
Permutation structures are obtained by taking automorphism groups
of first-order structures and we often regard a first-order
structure as a permutation structure without explicitly saying so.
Let $\pi:C\rightarrow W$ be a finite cover (Definition \ref{finite
cover}),
 we frequently use the notation $C(w)$ to denote the fibre
$\pi^{-1}(w)$ above $w$ in the cover $\pi:C\rightarrow W$.\\
We recall that the {\it{fibre group}} $F(w)$ of $\pi$ on $C(w)$ is
Aut$(C(w)/w)$, while the {\it{binding group}} $B(w)$ of $\pi$ on
$C(w)$ is Aut$(C(w)/W)$. It follows that the binding group is a
normal subgroup of the fibre group. If Aut$(W)$ acts transitively
on $W$, then all the fibre groups are isomorphic as permutation
groups, as are the binding groups. There is a continuous
epimorphism $\chi_w:\textrm{Aut}(W/w)\rightarrow F(w)/B(w)$ called
{\it{canonical epimorphism}} (Lemma 2.1.1 \cite{EMI}). Thus if $\Aut(W/w)$ has no proper
open subgroup of finite index, then
$F(w)=B(w)$.\\
Let $\pi_1:C_1\rightarrow W$ and  $\pi_2:C_2\rightarrow W$ be two finite covers of $W$. Then $\pi_1$ is said to be {\it{isomorphic}} to  $\pi_2$  if there exists a bijection $\alpha:C_1\rightarrow C_2$ with $\alpha(\pi_1^{-1}(w))=\pi_2(w)$ for all $w\in W$,  such that the induced map $f_\alpha:\Sym(C_1)\rightarrow \Sym(C_2)$  satisfies $f_\alpha(\Aut(C_1))=\Aut(C_2)$.\\

 The {\it{Cover
Problem}} is, given $W$ and data $(F(w),B(w),\chi_w)$, to
determine (up to isomorphism) the possible finite covers with these data.\\

If $C$ and $C'$ are permutation structures with the same domain
and $\pi:C\rightarrow W$,  $\pi':C'\rightarrow W$ are finite covers
with $\pi(c)=\pi'(c)$ for all $c\in C=C'$,  we say that $\pi'$ is
a {\it{covering expansion}} of $\pi$ if
Aut$(C')\leq$ Aut$(C)$.\\

Suppose that $C$ and $W$ are two permutation structures and
$\pi:C\rightarrow W$ is a finite cover. The cover is {\it{free}}
if
$$
\textrm{Aut}(C/W)=\prod_{w\in W}\textrm{Aut}(C(w)/W),
$$
that is, the kernel is the full direct product of the binding
groups.\\
The existence of a free finite cover with prescribed data depends
on the existence of a certain
continuous epimorphism.\\
Indeed, let $W$ be a transitive permutation structure and $w_0\in
W$. Given a permutation group $F$ on a finite set $X$, a normal
subgroup $B$ of $F$ and a continuous epimorphism
$$\chi: \textrm{Aut}(W/w_0)\rightarrow F/B,$$ then there exists a free finite cover
$\sigma:M\rightarrow W$ with fibre and binding groups at $w_0$
equal to $F$ and $B$, and such that the canonical epimorphism
$\chi_{w_0}$ is equal to $\chi$ . With these properties $\sigma$ is determined uniquely (see \cite{EMI}, Lemma 2.1.2).\\

A $\emph{principal}$ cover $\pi:C\rightarrow W$ is a free finite cover where the fibre and binding groups at each point are equal.
Free covers are useful in describing finite covers with given data
because every finite cover $\pi:C\rightarrow W$ is an expansion of
a free finite cover with the same fibre groups, binding groups and
canonical homomorphisms as in $\pi$ (see \cite{EMI}, Lemma
2.1.3).\\

Let's go back to Section 2. Using the language of finite covers, $\mathcal{F}$ is the set of the expansions of the principal finite covers of $\langle W, \Upsilon \rangle$, with all fibre groups and binding groups equal to  a given group $G$.\\ 
In the case when $G$ is a simple non-abelian regular group, our main theorem shows that the $\Upsilon$-congruences on $W$ describe (up to isomorphisms over $W$) the kernels of expansions belonging to $\mathcal{F}$.

\subsection{Bi-interpretability}

\begin{defn}
Two permutation structures are {\it{bi-interpretable}} if their
automorphism groups are isomorphic as topological groups.
\end{defn}
 For  a model-theoretic interpretation, if the permutation structures arise from $\aleph_0$-categorical
structures, see Ahlbrandt and Ziegler (\cite{AZ}). Usually
classification of structures is up to bi-interpretability.\\

Let $n\in\N$. Consider $\Omn$ as a first-order structure with automorphism group equal to $\Sym(\Om)$.

\begin{prop}\label{biinterpretability}Let $M_1:=\Delta\times \Omn$ and  $\pi_1:M_1\rightarrow \Omn$ be an expansion of a principal finite cover of $\Omn$ with all binding groups equal to a simple non-abelian finite group $G$ acting on $\Delta$. Let  $K_1$ be the kernel of $\pi_1$.\\
Suppose that the congruence classes which   $K_1$ determine have finite cardinality. Then, $\forall m>n$ there exists a permutation structure $M_2:=\Delta\times \Omm$ and a finite cover $\pi_2:M_2\rightarrow \Omega^{(m)}$  with  all fibre groups and binding groups equal to $G$ 
such that  $M_1$ is bi-interpretable with $M_2$
and the kernel $K_2$ of $\pi_2$ determines a $\Sym(\Om)$-congruence with equivalence classes of infinite cardinality .
\end{prop}
{\emph{Proof. }}
By the notation $M_1(\alpha)$, we mean the copy of $\Delta$ over the element $\alpha\in \Omn$. The kernel $K_1$, by Lemma \ref{pregeometry}, determines a $\Sym(\Om)$-congruence $\rho$ which, by hypothesis has equivalence classes of finite cardinality. Let $m$  be a positive integer greater then $n$ and $M_2$  be the set
$$
M_2=\{(w,m): w=(\alpha,c_1,\dots,c_{m-n}) \hbox{ and } m\in M_1(\alpha)\}
$$
where $\alpha \in \Omn$ and $c_1,\dots,c_{m-n}\in \Om\setminus \hbox{supp}(\alpha)$ and are all distinct. Obviously $M_2=\Delta\times \Omm$. Let $\mu_1:\Aut(M_1)\rightarrow \Sym(\Om)$ be the map induced by $\pi_1$  and $\Lambda$  be the subgroup of $\Sym(\Om)\times \Aut(M_1)$
$$
\Lambda=\{(g,\sigma) : g=\mu_1(\sigma)\}.
$$
Our claim is to show that $\langle M_2, \Lambda \rangle $ is a permutation structure and that $\pi_2:M_2\rightarrow \Omm$ given by $\pi_2(w,m)=w$ is a finite cover of $\Omm$ with  $F(w)=B(w)=G$ and kernel $K_2$ which determines a $\Sym(\Om)$-congruence with equivalence classes of infinite cardinality .\\
It is easy to check that $\Lambda$ is a permutation group on $M_2$ which preserves the partition  of $M_2$ given by the fibres of $\pi_2$.\\
 
We equip  $\Sym(\Om)\times \Aut(M_1)$ with the product topology. This topology coincides with the  topology of the pointwise convergence induced by   $\Sym(\Omm\times M_1)$ on $\Sym(\Om)\times \Aut(M_1)$.
The map $\Phi$ given by
$$\Sym(\Om)\times \Aut(M_1)\overset{p_1}{\rightarrow}\Sym(\Om)$$
  and the map $\Psi$ given by
$$\Sym(\Om)\times \Aut(M_1)\overset{p_2}{\rightarrow}\Aut(M_1)\overset{\mu_1}{\rightarrow}\Sym(\Om)$$
where $p_1$ and $p_2$ are the projections on the first and second component, respectively,   are continuous.   The permutation group $\Lambda$ is equal to the difference kernel
$$Z=\{(g,\sigma)\in\Sym(\Om)\times \Aut(M_1) : \Psi(g,\sigma)=\Phi(g,\sigma) \}$$
which, by Proposition 3 pag. 30  of \cite{Hig},  is closed in $\Sym(\Om)\times \Aut(M_1)$. Moreover, $\Sym(\Om)\times \Aut(M_1)$ is closed in $\Sym(\Omm\times M_1)$ and then  $\langle M_2, \Lambda \rangle $ is a permutation structure. The usual map induced by $\pi_2$
$$\mu_2: \Lambda \rightarrow \Sym(\Omm)$$ has image $\Sym(\Om)$. The kernel of $\mu_2$, which we denote by $K_2$, is
$$
K_2=\{(id,\sigma)\in \Lambda : \sigma\in K_1\}.
$$
Then $K_1\cong K_2$.
Let $(w,m)=(\alpha,c_1,\dots,c_{m-n}, m)\in M_2$ where $\alpha \in \Omn$ and $c_1,\dots,c_{m-n}\in \Om\setminus \hbox{supp}(\alpha)$ and are all distinct. Let $(id, \sigma)$ be an element in $\K_2$. If we restrict  it to the fibre over $w$, we see that it is the same as restricting $\sigma$ to the fibre over $\alpha$. Hence the binding group over $w$, $B_2(w)$, is clearly isomorphic to $G$. The  same holds for the fibre group:  let  $w=(\alpha,c_1,\dots,c_{m-n})$, then $F_2(w)$ over is the restriction of   the group $$\Aut(M_2/w)=\{(g,\sigma)\in \Lambda: g\in \Sym(\Om)_{((\alpha,c_1,\dots,c_{m-n})}\}$$  to the fibre over $w$. Since $g\in \Sym(\Om)_{((\alpha,c_1,\dots,c_{m-n})}\}$  then $g\in \Sym(\Om)_{(\alpha)}$. Hence $\sigma\in \Aut(M_1/\alpha)$ and so $F_2(w)$ is  isomorphic to $G$.

 Moreover, if we consider two points of $\Omm$, say $w=(\alpha, c_1,\dots,c_{m-n})$ and $w'=(\alpha', c'_1,\dots,c'_{m-n})$, with $\alpha\rho \alpha'$, we have that  $K_2(w,w')\cong G$. Vice versa if $K_2(w,w')\cong G$, it means that  $K_1(\alpha, \alpha')\cong G$. Then the $\Sym(\Om)$-congruence, $\rho'$, that $K_2$ determines is given by $w\rho' w'$ if and only if $\alpha\rho \alpha'$. In the equivalence class of $w=(\alpha, c_1,\dots,c_{m-n})$ for instance there are all the elements  of the form $(\alpha, c'_1,\dots,c'_{m-n})$,  with $c_1,\dots,c_{m-n}\in \Om\setminus \hbox{supp}(\alpha)$ and  pairwise distinct. Then the equivalence classes of $\rho'$ are of infinite cardinality.\\
 Next we check the bi-interpretability. We consider the map
 $$
 \begin{array}{cccc}
  \beta:& \Lambda &\rightarrow & \Aut(M_1)\\
   & (g,\sigma)& \mapsto & \sigma
  \end{array}
 $$
The kernel of $\beta$ is ker$\beta=\{(g,id)\in \Lambda : g=\mu_1(id)\}$.  Then $\beta$ is injective. It is also surjective since, given $\sigma\in \Aut(M_1)$, $(\mu_1(\sigma), \sigma)\in \Lambda$. Clearly  the inverse map is given by $\beta^{-1}(\sigma)=(\mu_1(\sigma), \sigma)$.\\ 
 It is a topological isomorphism. Indeed, take a basic open neighbourhood of the identity  in $\Aut(M_1)$, say $\Aut(M_1)_{(\Gamma)}$, where $\Gamma=\{m_i\}_{i\in I}$ is a finite set of $M_1$. Each $m_i\in M_1(\alpha_i)$. Then $$\beta^{-1}(\Aut(M_1)_{(\Gamma)})=\{(\mu_1(\sigma), \sigma) : \sigma\in \Aut(M_1)_{(\Gamma)}\}.$$
 For each $\alpha_i$, we choose $c^i_1,\dots,c^i_{m-n}\in \Om$  such that  $w_i=(\alpha_i,c^i_1,\dots,c^i_{m-n})$ is a fulfillment of $\alpha_i$ to an element of $\Omm$. The map 
 $$
 \begin{array}{cccc}
 \beta^{-1}: & \Aut(M_1) & \rightarrow & \Sym(\Om)\times \Aut(M_1)\\
 & \sigma & \mapsto & ( \mu_1(\sigma), \sigma)
 \end{array}
 $$
  is continuous. The image of $\beta^{-1}$ is $\Lambda$ and $\Lambda$ has the topology induced by  $\Sym(\Om)\times \Aut(M_1)$, then $ \beta^{-1}:\Aut(M_1)\rightarrow  \Lambda$ is continuous. Hence, we have proved the bi-interpretability.\cvd

\subsection{Almost-free finite covers}
Let $W$  be a transitive structure and  $\rho$ be  an $\Aut(W)$-congruence on $W$. Given a $\rho$-equivalence class $[w]$, we denote by $C([w])=\bigcup_{w_i\in [w] }C(w_i)$, by $F([w])$ the permutation group induced by $\Aut(C/\{[w]\})$ on  $C([w])$, and by $B([w])$ the permutation group induced by the kernel of $\pi$ on $C([w])$.\\ 
Note that $B([w])\unlhd F([w])$.

\begin{lemma}\label{chi-map}
Suppose that $W$ is a transitive structure and $\rho$  an $\Aut(W)$-congruence on $W$. Let  $\pi:C\rightarrow W$ be a finite cover. Then, for every $\rho$-class $[w]$ in $W$
\begin{enumerate}
\item there exists a finite-to-one surjection 
$$
\pi_{[w]}:C([w])\rightarrow [w]
$$
such that its fibres form an $F([w])$-invariant partition of $C([w])$;
\item there is a continuous epimorphism
$$
\chi_{[w]}:\Aut(W/\{[w]\})\rightarrow F([w])/B([w]).
$$
\end{enumerate}

 \end{lemma}

\emph{Proof.} The first point is clear.\\   
The second point require  a little proof. Let $g\in \Aut(W/\{[w]\})$. Then there exists $h\in \Aut(C/\{[w]\})$ which extends $g$.   Let $\psi:\Aut(W/\{[w]\})\rightarrow
\Aut(C/\{[w]\})/\Aut(C/W)$ be the map defined by
$\psi(g)=h\Aut(C/W)$. This map is well defined. Suppose that also $\bar{h}$ extends $g$. Then $h^{-1} \bar{h}\in \Aut(C/W)$ and so $h\Aut(C/W)=\bar{h}\Aut(C/W)$.  Consider the restriction to the set of fibres over $\{[w]\}$.  So we have
a map 
$\xi_{[w]}:\Aut(C/\{[w]\})/ \Aut(C/W)\rightarrow
\Sym(C([w])/B([w])$, given by $\xi_{[w]}(h\Aut(C/W))=h_{|C([w])}B([w])$, which is clearly onto on $F([w])/B([w])$. Let $g\in \Aut(W/\{[w]\})$.
We define $\chi_{[w]}(g):=\xi_{[w]}\psi(g)$. 
In order to prove that $\chi_{[w]}$ is continuous, we show that
$\psi$ and $\xi_{[w]}$ are continuous.\\
The restriction map  $\xi_{[w]}$ is continuous by Lemma 1.4.1 of \cite{EMI}. Consider
$\Sym(C([w])$ with the topology of pointwise
convergence and $\Sym(C([w])/B([w])$ with the quotient
topology.  Let $\mu_{| \Aut(C/\{[w]\})}: \Aut(C/\{[w]\})\rightarrow \Aut(W/\{[w]\})$ be the map induced by $\mu$. Since $[w]$ is a $\rho$-equivalence class $ \Aut(C/\{[w]\})$ is an open subgroup of $\Aut(C)$. Indeed, let $c\in C([w])$. Take $h\in \Aut(C/c)$. Then $h( C([w]))= C([w])$. If $g=\mu(h)$, we have $g(w)=w$, and being $[w]$ a $\Aut(W)$-congruence class, this implies that $g([w])=[w]$. Hence $h\in \Aut(C/\{[w]\})$. By Lemma \ref{trivialita} we have that $ \Aut(C/\{[w]\})$ is an open subgroup of $\Aut(C)$. By the same reasoning we get that   $\Aut(W/\{[w]\})$ is open in $\Aut(W)$.  Now, since $\mu$ is open also $\mu_{| \Aut(C/\{[w]\})}$ will be open. Hence by Proposition 1, pag 21 of \cite{Hig}, we have the continuity of $\psi$.
\cvd

\begin{defn}

Let $\pi:C\rightarrow W$ be a finite cover of $W$, $w\in W$,  with binding groups isomorphic to a group  $G$ and kernel $K$. We shall say that $\pi$ is   {\bf{almost free} } with respect to $\rho$
  if

\begin{enumerate}
\item $K([w])\cong G$ for each $[w]\in W/\rho$ 
\item  $K(w_1,w_2)\cong G\times G$ for each $w_2\notin [w_1].$
\end{enumerate}
\end{defn}
A class of almost free finite cover is the set of the expansions of the free finite covers with  binding groups isomorphic to a simple non-abelian group $G$.\\

Let $R:=W/\rho$. Given a transitive structure $W$ and an $\Aut(W)$-congruence $\rho$, naturally we have an induced map 
$$
M:\Aut(W)\rightarrow \Sym(R).
$$
The map $M$ is continuous, but the image of $\Aut(W)$ by $M$  is not necessarily closed. The following counterexample is due to Peter Cameron (private communication).\\

Take the generic bipartite graph B, and consider the group G of
automorphisms fixing the two  bipartite blocks, acting on the set of
edges of the graph. On the set of edges there are two equivalence
relations, "same vertex in the first bipartite block", and "same vertex
in the second bipartite block". Clearly G is precisely the group
preserving these two equivalence relations, and so is closed. But the
group induced on the set of equivalence classes of each relation is
highly transitive and not the symmetric group, therefore not closed.\\

 \begin{prop}\label{teoremino} 
Let $W$ be a transitive structure and $\rho$ an $\Aut(W)$-congruence on $W$. We suppose that the following assumptions hold:

\begin{enumerate}
 
\item Let $F$ be a closed permutation group on a   set $X$. Fix $w_0\in W$ and let $[w_0]$ be the $\rho$-equivalence class of $w_0$.\\
 Suppose that there exists a finite -to-one  surjection $$\sigma:X\rightarrow [w_0]$$  such that the fibres form an $F$-invariant partition of $X$ and that the induced map $T: F\rightarrow  \Sym([w_0])$ has image $\Aut(W/\{[w_0]\}_{|\{[w_0]\}}$. Let $B$ the kernel of $T$.
 \item The map $T$ induces a map
$$
\chi:\Aut(W/\{[w_0]\})\rightarrow F/B
$$
defined as $\chi(g)=hB$, where $h\in F$ and $T(h)=g_{|[w]}$. Assume that  $\chi$ is continuous.

\item Let $G$ be the permutation group induced by $B$ on $\sigma^{-1}(w_0)$ . Suppose that $B$ is isomorphic to $G$.
\item Assume that the map $M$ is injective,  open and with closed image.

\end{enumerate}

Then there exists an almost free finite cover $\pi_0$ of $W$ with respect to $\rho$  with binding groups isomorphic to $G$, $F([w_0])=F$, $B([w_0])=B$ and map $\chi_{[w_0]}$ equal to $\chi$. Moreover, if  $\tilde{\pi}_0$ is an almost free finite cover with respect to $\rho$ with $F([w_0])$  and $B([w_0])$ isomorphic as permutation groups to $F$ and $B$ respectively, and $\chi_{[w_0]}$ equal to $\chi$ (up to isomorphism), then  $\tilde{\pi}_0$ is isomorphic over $W$ to $\pi_0$. 
\end{prop}
\emph{Proof.} This is an  application of  Lemma 2.1.2 in \cite{EMI}. In this proof we will deal with a map with all the properties of a finite cover but  the finitess condition on the fibres (hence we allow the cover to having fibres of infinite cardianality). We shall call such a map a cover.\\ 
 We give to $R$ the first-order structure with automorphism group  the image of $M$. Let $r_0=[w_0]$. We have that $$M^{-1}:\Aut(R/r_0)\rightarrow \Aut(W/\{[w_0]\})$$ is continuous. Then we have a continuous map  $\chi:\Aut(R/r_0)\rightarrow F/B$.\\
Since we are going to use a slightly changed version  of the proof of Lemma 2.1.2 in \cite{EMI} and then to use specific steps out of it, we are going to give the general lines of the proof for the use of the reader. For the details we  refer  to the book \cite{EMI}.\\ 

We are going to sketch the proof of the following statement: Let $R$ be a transitive permutation structure and $r_0\in R$ .  Let $F$ be a closed permutation group on a  set $X$, and $B$ be a normal subgroup of $F$. Suppose there is a continuous epimorphism $\chi:\Aut(R/r_0)\rightarrow F/B$. Then there exists a cover $\pi:M\rightarrow R$ with fibre group  and binding group at $r_0$ respectively equal to $F$ and $B$ and canonical epimorphism at $r_0$ equal to $\chi$. Moreover, if $\nu$ is a free  cover with with $F(r_0)$  and $B(r_0)$ isomorphic as permutation groups to $F$ and $B$ respectively, and $\chi_{r_0}$ equal to $\chi$ (up to isomorphism), then  $\nu$ is isomorphic over $R$ to $\pi$.\\  

The proof is made essentially in three steps. First  the following cover  is constructed.\\
Let $C$ be the set of left cosets of $\Ker\chi$ in
$\Aut(R)$. Consider the map $\theta:C\rightarrow R$ given by
$\theta(g\Ker\chi)=gr_0$. The permutation group $\Aut(R)$
induces a group of permutation on $C$. The induced group is a closed subgroup of $\Sym(C)$ and so we can consider $C$
as a relational structure with automorphism group isomorphic to
$\Aut(R)$. Then the map $\theta$ is a cover
 with trivial kernel.\\
 Let $Y=\theta^{-1}(r_0)\cup X$. Put on $Y$ the relational
 structure which has as automorphism group $F$: the action of $h\in F$ on
 $m\in \theta^{-1}(r_0)$ is $h(m)=(\chi^{-1}(hG))(m)$. For
 every $r\in R$ choose $g_r\in \Aut(R)$ such that $g_rr=r_0$
 (with $g_{r_0}=id$). Then $g_r(\theta^{-1}(r))=
 \theta^{-1}(r_0)$ and it induces an embedding
 $\eta_r:\theta^{-1}(r)\rightarrow Y$.\\

The second step is the following: we built a  cover $\pi':M'\rightarrow R$, where the domain of
 $M'$ is made of the disjoint union of $R$, $C$ and $R\times
 Y$ and $\pi'$ is the identity on $R$, $\theta$ on $C$, the
 projection to the first coordinate on $R\times Y$. We also have an
 injection $\tau:C\rightarrow R\times Y$ given by
 $\tau(c)=(r,\eta_r(c))$, whenever $\theta(c)=r$. Moreover, the
 structure of $M'$ is made up of the original structure on $R$ and
 $C$ and for each $n$-ary relation $R$ on $Y$ we have an
 $n$-ary relation $R'$ on $R\times Y$ given by
 $$
R'((r_1,y_1),\cdots,(r_n,y_n))\hbox{ iff }r_i=r_j,\, \forall i,j
\hbox{ and } R(y_1,\cdots,y_n).
 $$
Now we see how to extend
an automorphism of $R$  to a permutation of $M'$ which preserves the above structure.\\
Let $g\in \Aut(R)$, then we get an automorphism of $C$. Let $g r_1=r_2$, then  via $\tau$ we have a
bijection from $\{r_1\}\times \theta^{-1}(r_0)$ to $\{r_2\}\times
\theta^{-1}(r_0)$. In fact, let $\bar{g}\Ker\chi\in \theta^{-1}(r_0)$,
then $\tau g \tau^{-1}(\bar{g}\Ker\chi)=g_{r_2}gg^{-1}_{r_1}\bar{g}\Ker\chi$. Since
$g_{r_2}gg^{-1}_{r_1}\in \Aut(W/\{r_0\})$, if we choose a
representative $h$ in the class $\chi(g_{r_2}gg^{-1}_{r_1})$ then
$h(\bar{g}\Ker\chi)=\tau g\tau^{-1}(\bar{g}\Ker\chi)$ and this extends to a permutation
$\beta(r,g)$ of $Y$. If we also denote by $\beta(r,g)$ the induced
map from $r\times Y$ to $gr\times Y$, then
$\omega(g)=g\cup \bigcup_{r\in R}\beta(r,g)$ is a permutation of $M'$
which preserves the structure we put above on $M'$ and extends $g$.\\
Let $\pi$ the restriction of $\pi'$ to
$M=R\times X$ considered as permutation structure with
$\Aut(M')$ acting. Then 
$\pi:M\rightarrow R$ is a free cover of $R$ and kernel isomorphic to $G^R$.\\
Now the uniqueness, the third step.  Let $\gamma:N(w_0)\rightarrow X$ be the bijection which gives rise to the isomorphism (we call it  $\tilde{\gamma}$) as permutation groups between $F$  ($B$)  and $F(w_0)$ ($B(w_0)$).  Let $\nu:N\rightarrow R$ be a 
cover with $F(w_0)$ and $B(w_0)$ isomorphic as permutation groups to $F$ and $B$ respectively and $\chi_{w_0}=\tilde{\gamma}\circ\chi$. For each $r\in R$, $g_{r}$
can be extended to an automorphism $\hat{g}_{r}\in \Aut(N)$.
 We define the map $\beta:N \rightarrow
R\times X$ in the following way: if $n\in \nu^{-1}(r)$, 
define $\beta(n):=(r,\gamma(\hat{g}_{r}(n)))\in R\times X$. As it is shown in  Lemma 2.1.2 in \cite{EMI}, this is a bijection which gives rise to an isomorphism of covers.\\

Let $g_r\in \Aut(R)$ be the permutations used above for constructing the free  cover $M$.  Then we construct a finite cover of $W$ in the following way. Consider the set 

$$C_0:=\{(w, k)\, :\, w\in r \hbox{  and }  k\in \sigma^{-1}(M^{-1}(g_r)(w))\}$$
The map $\pi_0:C_0\rightarrow W$ given by $\pi_0(w,k)=w$ is a finite-to-one surjection.\\
Let $\alpha: R\times X\rightarrow C_0$  be the map defined in the following way: let $k\in X$, then there exists $w\in [w_0]$ such that  $k\in \sigma^{-1}(w)$. We define  $$\alpha(r, k):=(M^{-1}(g^{-1}_r) w, k).$$ This is a bijection. Let  $f_\alpha: \Sym(M)\rightarrow \Sym(C_0)$ be the induced map by $\alpha$. The image by $f_\alpha$ of $\Aut(M)$ is closed in $\Sym(C_0)$. We denote it by  $\Aut(C_0)$.\\  
Let $C_0(w)$ be the fibre over $w$ of $\pi_0$. If $w\in r_1$ then $C_0(w)=\sigma^{-1}(M^{-1}(g_{r_1})( w)$. We have that $\alpha^{-1}C_0(w)=(r_1,\sigma^{-1}(M^{-1}(g_{r_1}) (w))$.\\ 
Take an element $g$ of $\Aut(M)$. We are going to show that $\alpha g \alpha^{-1}$ preserves the partition of $C_0$ given by the fibres of $\pi_0$.\\ 
Let $\bar{g}\in \Aut(W)$ such that $M(\bar{g})$ is the induced permutation on $R$ by $g$. If $M(\bar{g})r_1=r_2$, there exists $f\in F$   such that $$g(r_1,\sigma^{-1}(M^{-1}(g_{r_1}) w))=f(\sigma^{-1}(M^{-1}(g_{r_1}) w)=\sigma^{-1}(M^{-1}(g_{r_2})\bar{g} w).$$ 
By the proof of Lemma 2.1.2 in \cite{EMI}, we see that the element $f$ is a representative of the class  $\chi(M^{-1}(g_{r_2})\bar{g}M ^{-1}(g^{-1}_{r_1}))$. Hence $g(r_1,\sigma^{-1}(M^{-1}(g_{r_1}) w))=(r_2, \sigma^{-1}(M^{-1}(g_{r_2})\bar{g} w))$ and then
$$
\alpha g\alpha^{-1}C_0(w)=C_0(\bar{g} w),
$$
i.e. the fibres of $\pi_0$ form an $\Aut(C_0)$-invariant partition of $C_0$.\\
Let  $\mu_0:\Aut(C_0)\rightarrow \Sym(W)$ be the induced homomorphism. Take an element $g\in \Aut(W)$ and an extension $\tilde{g}\in \Aut(M)$ of $M(g)$. The argument above shows as well  that the Im$\mu_0$ is equal to $\Aut(W)$. The kernel of $\mu_0$ is $\alpha \Ker\pi\alpha^{-1}$.  It is isomorphic to $G^R$. Since $\Ker\pi$ induced on $\sigma^{-1}(w)$ and on $X$ is isomorphic to $G$,  then $\Ker\pi_0$ induced on any fibre of $\pi_0$ and on $C_0([w_0])$ is isomorphic to $G$ as well.  So we have an almost free finite cover  $\pi_0:C_0\rightarrow W$ as required.\\

Let $\nu_0:N_0\rightarrow W$ be a finite cover of $W$ with binding groups isomorphic to a  finite group $G$ with kernel isomorphic to $G^R$ and with $B([w_0])\cong B$ and $F([w_0])\cong F$  as permutation groups. Suppose that $\chi_{[w_0]}$ is equal to $\chi$.\\
Let $\Delta[w]:=N_0([w])$ and $\Delta=\cup_{[w]\in R} \Delta[w]$. Let $\nu:\Delta\rightarrow R$ given in the obvious way by $\nu(\delta)=r$ if $\delta \in \Delta[w]$ and $[w]=r$. The group $\Aut(N_0)$ acts on it and can be taken as automorphism group of $\Delta$.\\
The fibres of $\nu$ form a partition of $\Delta$ invariant under the action of $\Aut(N_0)$. Indeed,  let $g\in \Aut(R)$, consider $M^{-1}(g)$ which extends to $\bar{g}\in \Aut(N_0)$. Then, if $\delta\in \Delta[w]$ there exists $n\in [w]$ such that $\delta\in N_0(n)$ and  $\bar{g}\delta\in N_0 (M^{-1}gw)\subseteq\Delta(g[w])$.\\  
The fibre group at $r_0$ is equal to $F[w_0]$ and the binding group at $r_0$ is equal to $B([w_0])$. The map $\chi_{r_0} :\Aut(R/r_0)\rightarrow F([w_0])/B([w_0])$ is exactly the composition of $M^{-1}:\Aut(R/r_0)\rightarrow \Aut(W/\{[w_0]\})$ and $\chi_{[w_0]}$. Since the data of $\pi$ and $\nu$ are the same up to isomorphism, by  Lemma 2.1.2 in \cite{EMI}  $\nu$ and $\pi$ are isomorphic over $R$ via the bijection $\beta(\delta)=([w], \gamma(\hat{g_r}(\delta)))$, if $\delta\in \Delta[w]$ and $\hat{g_r}\in \Aut(N_0)$ is an extension of $M^{-1}g_r$.\\ 
Let $\delta \in N_0(w)$ (so $\hat{g_r}\delta\in N_0(M^{-1}(g_r)w)$). Consider the bijection
$$
\begin{array}{ccccccc}
  N_0 &\overset{id}{\rightarrow} & \Delta  & \overset{\beta}{\rightarrow} & M & \overset{\alpha}{\rightarrow} &  C_0  \\
 \delta & \mapsto  &  \delta  & \mapsto & ( [w], \gamma(\hat{g_r}\delta)) &\mapsto &  (w, \gamma(\hat{g_r}\delta))
\end{array}
$$
Then $\alpha\beta\Aut(N_0)\beta^{-1}\alpha^{-1}=\alpha\Aut(M)\alpha^{-1}=\Aut(C_0)$, i.e.  $\Aut(N_0)$ and $\Aut(C_0)$ are isomorphic over $W$. \cvd  

\begin{ex}
 Let $W$ be a transitive structure, $w_0\in W$, and $\rho$ be an  $\Aut(W)$-congruence on $W$. Assume that the permutation group induced by $\Aut(W/\{[w_0]\})$ on $[w_0]$, which we shall denote by $A$, is closed in $\Sym([w_0])$. Moreover suppose that the map $M$ is injective,  open and with closed image, as in Proposition \ref{teoremino}. Let $G$ be a  finite permutation group acting on a set $L$ . There always exists an almost-free finite cover. In order to see it, consider the wreath product $GWr_{[w_0]}A$ acting  in the usual way on $[w_0]\times L$.\\ 
Let $\sigma: [w_0]\times L\rightarrow [w_0]$ given by $\sigma(w,l)=w$.  Denote by $B_1$ the diagonal subgroup of $G^{[w_0]}$: it is normalized by $A$ and so we can make the semidirect product $F_1:=B_1 \rtimes A$. This is closed by Proposition \ref{closedcompact}. Using the notation of Proposition \ref{teoremino} we have that $\chi$ is the homomorphism induced by restriction on $[w_0]$. Since $\chi$  is continuous, the hypothesis of Proposition \ref{teoremino} are satisfied and so we have an almost free finite cover $\pi:W\times L\rightarrow W$. We note that the automorphism group $\Aut(W\times L)$,  which we have got, is equal to $K_\rho \rtimes \Aut(W)$ (using the notation of Theorem \ref{congruenze e Kernel}).\\
 
Now suppose that  $G$ is a simple non-abelian finite permutation group acting on itself by conjugation (so $G=L$). Next we give an example of an almost free finite cover with respect to $\rho$, not isomorphic to $\pi$,   with kernel equal to $\Ker\pi$.\\

Let $\pi:W\times G\rightarrow W$ be the cover that we have built above. Using the topological results in section 1.4 of \cite{EMI} we have that the map $T:F_1\rightarrow A$ is continuous, maps closed subgroups to closed subgroups and it is open. Then the isomorphism map $S:A\rightarrow F_1/B_1$ is  a topological isomorphism.\\
Since $B([w_0])=B_1\cong G$,  by conjugation of $G$ by elements of $F([w_0])=F_1$ we get a map $\gamma:F([w_0])/G\rightarrow \Out(G)$. The image of $\gamma$ is $H/G$, for some $H\leq \Aut(G)$. Composing  $S$ with $\gamma$, we have a map $$\bar{S}:A\rightarrow H/G.$$\\
In order to prove that $\gamma$ is continuous we have to show that the kernel of $\gamma$ is open. The kernel of $\gamma$ is $C(G)_{F([w_0])}G/G$, where $C(G)_{F([w_0])}$ is the centralizer of $G$ in $F([w_0])$. The group $G$ is finite and hence closed in $F([w_0])$. Its orbits on $[w_0]\times G$ are finite and so it is also compact. Moreover, $C(G)_{F([w_0])}$  is closed. \\ 
By Proposition \ref{closedcompact} we have that $C(G)_{F([w_0])}G$ is closed in $F([w_0])$. Since it has finite index in $F([w_0])$, $C(G)_{F([w_0])}G$ is open $F([w_0])$ and hence $C(G)_{F([w_0])}G/G$ is open in $F([w_0])/G$.\\
Let  $P:H\rightarrow H/G$ be the quotient map and 

$$
F_2:= \{(\sigma, h): \sigma \in A, \, h\in H \hbox{ and } P(h)=S(\sigma)\}
$$ 
be the fibre product between $A$ and $H$. This is a permutation group on $[w_0]\times G$ with action given by: $(\sigma, h)(w,g)=(\sigma w, h(g))$. By the same reasoning as in Proposition \ref{biinterpretability}, we have that $F_2$ is closed in $\Sym([w_0]\times G)$.\\ 
 
The group $B_2:=\{(id, g): id \in \Sym(\{[w_0]\}), \, g\in G\}$ is a normal subgroup of it. Let 
$\chi:Aut(W/\{[w_0]\})\twoheadrightarrow F_2/B_2$ be the map given by
 $$\chi (g)=(g_{|[w_0]},h)B_2,$$ where $h$ belongs to the coset $S(g_{|[w_0]})$.  The map $\chi$ is well defined. Moreover, $\chi$ is continuous, since $S$ is continuous.\\
 Let $$\sigma:[w_0]\times G\rightarrow [w_0]$$ be the projection on the first component.  The induced map $F_2\rightarrow \Sym([w_0])$ has image $A$. Hence, by Proposition \ref{teoremino}, we can build an almost-free finite cover $\pi_0$ w.r.t $\rho$ with binding groups isomorphic to $G$. Note that the kernel is equal to $K_\rho$. \end{ex}

\subsection{Problems}
We described in an explicit way the kernels of expansions of the
free finite cover of $\langle\Omn,\, \Sym(\Om)\rangle$, when the
fibre groups and the binding groups are both isomorphic to a
simple non-abelian finite group $G$.\\

1. What happens for finite covers where the base structure is a
Grassmannian of a vector space over a finite field?\\
2. What happens for finite covers of $\Omn$ if the fibre groups
and the binding groups are isomorphic to a simple abelian group?
Here one would need to work with the closed $\Sym(\Om)$-submodules
of $\mathbb{F}_p^{\Omn}$. We remind that the case where the base
permutation structure is $\langle\OmN, \, \Sym(\Om) \rangle$ was
solved by Gray (\cite{Gr}).

\vspace{1cm}

\begin{center}
\begin{large}
ACKNOWLEDGMENTS
\end{large}
\end{center}
The author wishes to express her thanks to  D. M. Evans for
several stimulating conversations and hospitality at UEA  and to
O. Puglisi for many helpful
 suggestions related to this paper.

\end{document}